\documentclass[12pt]{article}
\usepackage{amssymb}

\newtheorem{theorem}{Theorem}
\newtheorem{proposition}{Proposition}
\newtheorem{lemma}{Lemma}

\newcommand{\be}{\begin{equation}}
\newcommand{\ee}{\end{equation}}

\newcommand{\bea}{\begin{eqnarray}}
\newcommand{\eea}{\end{eqnarray}}

\newcommand{\hst}[1]{\rule{#1}{0mm}}
\newcommand{\ty}{\hspace{0.04em}}
\newcommand{\tty}{\hspace{0.06em}}

\newcommand{\sect}[1]{\setcounter{equation}{0}\section{#1}}

\newcommand{\reff}[1]{(\ref{#1})}

\newcommand{\dis}{\displaystyle}

\newcommand{\De}{\Delta}

\newcommand{\al}{\alpha}
\newcommand{\bet}{\beta}

\newcommand{\de}{\delta}

\newcommand{\ve}{\varepsilon}

\newcommand{\lam}{\lambda}

\newcommand{\vp}{\varphi}
\newcommand{\om}{\omega}

\newcommand{\la}{\langle}
\newcommand{\ra}{\rangle}
\newcommand{\rar}{\rightarrow}
\newcommand{\lra}{\longrightarrow}
\newcommand{\ti}{\times}
\newcommand{\op}{\oplus}

\newcommand{\ot}{\otimes}

\newcommand{\KK}{\mathbb K}

\newcommand{\ZZ}{\mathbb Z}

\newcommand{\CC}{\mathbb C}

\newcommand{\fg}{{\mathfrak g}}
\newcommand{\fgz}{{\mathfrak g}_{\overline0}}
\newcommand{\fgo}{{\mathfrak g}_{\overline1}}

\newcommand{\fs}{{\mathfrak s}}

\newcommand{\fz}{{\mathfrak z}}

\newcommand{\cI}{{\cal I}}
\newcommand{\cP}{{\cal P}}
\newcommand{\A}{{\cal A}}
\newcommand{\Ac}{{\cal A}^{\circ}}
\newcommand{\Ae}{{\cal A}_e}

\newcommand{\U}{{\cal U}}
\newcommand{\Uc}{{\cal U}^{\circ}}
\newcommand{\Ue}{{\cal U}_e}
\newcommand{\Uec}{{\cal U}_e^{\circ}}
\newcommand{\cO}{{\cal O}}
\newcommand{\B}{{\cal B}}
\newcommand{\cS}{{\cal S}}
\newcommand{\ol}{\overline}
\newcommand{\Oo}{\overline1}
\newcommand{\Oz}{\overline0}

\newcommand{\Ze}{\{0\}}
\newcommand{\cd}{\cdot}

\newcommand{\comp}{\!\stackrel{\textstyle{\hst{0ex} \atop \circ}}{\hst{0ex}}\!}
\newcommand{\id}{{\rm id}}

\newcommand{\tint}{{\textstyle\int}}
\newcommand{\tinte}{{\textstyle{\int_e}}}
\newcommand{\Hom}{{\rm Hom}}

\begin{document}

\begin{titlepage}

\newlength{\ppn}
\newcommand{\defboxn}[1]{\settowidth{\ppn}{#1}}
\defboxn{BONN--TH--2000--10}

%\mbox{}
\hspace*{\fill} \parbox{\ppn}{BONN--TH--2000--10 \\
                              November 2000        }

%\vspace{10mm}
\vspace{18mm}

\begin{center}
{\LARGE\bf Integration on Lie supergroups} \\[2.5ex]
{\large\bf A Hopf superalgebra approach} \\
\vspace{15mm}
{\large M. Scheunert} \\
\vspace{1mm}
Physikalisches Institut der Universit\"{a}t Bonn \\
Nu{\ss}allee 12, D--53115 Bonn, Germany \\
\vspace{5mm}
{\large R.B. Zhang} \\
\vspace{1mm}
School of Mathematics and Statistics \\
University of Sydney,
Sydney, NSW 2006, Australia \\
\end{center}

\vspace{15mm}
\begin{abstract}
\noindent
For a large class of finite-dimensional Lie superalgebras (including the
classical simple ones) a Lie supergroup associated to the algebra is defined
by fixing the Hopf superalgebra of functions on the supergroup. Then it is
shown that on this Hopf superalgebra there exists a non-zero left integral.
According to a recent work by the authors, this integral is unique up to
scalar multiples.
\end{abstract}

\vspace{\fill}
\noindent
math.RA/0012052

\end{titlepage}

\setcounter{page}{2}

\small\normalsize

\renewcommand{\thesection}{\arabic{section}.\hspace{-0.5em}}
\sect{Introduction \vspace{-1ex}}
\renewcommand{\thesection}{\arabic{section}}
The present work is a direct sequel to a recent paper by the authors
\cite{SZh}, in which we have studied the integrals on Hopf superalgebras,
with special consideration of the Hopf superalgebras of functions on
quantum (i.e., $q$-deformed) or classical (i.e., undeformed) Lie
supergroups. The reasons why we got interested in this topic and a list
of pertinent references can be found in that paper.

The uniqueness (up to scalar multiples) of integrals on Hopf superalgebras
could be proved by reducing the problem to the non-graded case. As expected,
to prove the existence of non-trivial integrals is a more difficult task.
Unfortunately, in Ref.~\cite{SZh}, we were only able to prove the existence
for some special classes of Hopf superalgebras. In particular, for most of
the quantum or classical orthosymplectic supergroups, the problem could not
be solved. Now, in the quantum case, there are clear indications that our
approach does not work for these examples (at least not in the form presented
in Ref.~\cite{SZh}). However, in the undeformed case, there are no such
reasons, and we were pretty sure that non-trivial integrals should exist at
least for the Lie supergroups associated to the basic classical simple Lie
superalgebras. The main object of the present work is to show that this is
indeed the case, in fact, we can prove this for even more general algebras, 
including all classical simple Lie superalgebras.

The present paper is set up as follows. In Section 2 we remind the reader
of the basic definitions and of the uniqueness result obtained in
Ref.~\cite{SZh}. Section 3 contains some comments on the finite dual of a
Hopf superalgebra that will be useful later on. In Section 4 we recall our
approach to the construction of integrals on Lie supergroups, and we also
explain how the Hopf superalgebra of functions on a Lie supergroup will be
chosen. The basic existence theorem will be proved in Section 5. A brief
discussion in Section 6 closes the main body of the paper.

As the reader will see, the theory of so-called Frobenius extensions of
rings \cite{NTs, BFa} will play a crucial role. In view of the importance of
our existence result for our work, we think it might be useful to include
some results of that theory. This will be done in the Appendix.

We close this introduction by explaining some of our conventions. In the
main body of the paper, the base field will be the field $\CC$ of complex
numbers (although an arbitrary algebraically closed field of characteristic
zero would do as well). Our notation for $\ZZ_2\tty$-graded algebraic
structures is fairly standard. The two elements of $\ZZ_2$ are denoted by
$\Oz$ and $\Oo$. Unless stated otherwise, all gradations considered in this
work will be $\ZZ_2\tty$-gradations. For any superspace, i.e.
$\ZZ_2\tty$-graded vector space $V = V_{\Oz} \op V_{\Oo}$\,, we define the
gradation index $[\;\,]:\, V_{\Oz} \cup V_{\Oo} \rar \ZZ_2$ by $[x]=\alpha$
if $x \in V_{\alpha}\tty$, where $\alpha \in \ZZ_2\,$. Unless stated
otherwise, {\em all algebraic notions and constructions are to be understood
in the {\em super sense}, i.e., they are assumed to be consistent with the
$\ZZ_2\tty$-gradations and to include the appropriate sign factors.}
\newpage

\renewcommand{\thesection}{\arabic{section}.\hspace{-0.5em}}
\sect{Definition and uniqueness of integrals \vspace{-1ex}}
\renewcommand{\thesection}{\arabic{section}}
Let $\A$ be a Hopf superalgebra with comultiplication $\Delta\tty$, counit
$\ve$, and antipode $S$. A {\em left integral} $\int^\ell$ on $\A$ is an
element of $\A^*$ (the dual of the graded vector space $\A$), such that
\begin{equation}
  (\id_{\A}\,\ot\tint^\ell) \comp \Delta
                                   \,=\, 1_{\A}\tint^\ell. \label{intleft}
\end{equation}
Equivalently, this means that
\begin{equation}
  a^{\ast} \cd \tint^\ell =\, a^{\ast}(1_{\A}) \tint^\ell \quad
                    \mbox{for all $a^{\ast} \in A^{\ast}$}, \label{intexpl}
\end{equation}
where the dot denotes the multiplication in $\A^{\ast}$ deduced from
$\Delta$\,. A right integral $\int^r\in\A^*$ on $\A$ is defined by
a similar requirement
\begin{equation}
  (\tint^r \otimes\,\id_{\A}) \comp \Delta \,=\, 1_{\A}\tint^r.
\end{equation}

Let $\A^{\rm aop,cop}$ be the Hopf superalgebra opposite to $\A$ both
regarded as an algebra and a coalgebra. Then a linear form
$\int \in \A^{\ast}$ is a left/right integral on $\A$ if and only if it
is a right/left integral on $\A^{\rm aop,cop}$. In particular, if the
antipode $S$ of $\A$ (and hence of $\A^{\rm aop,cop}$) is invertible,
then $S^{\pm 1}$ are isomorphisms of $\A$ onto $\A^{\rm aop,cop}$, and
hence $\int$ is a left integral on $\A$ if and only if
$\int \comp S^{\pm 1}$ are right integrals on $\A$\,. Thus we only need to
consider left integrals (or right integrals).

In Ref.~\cite{SZh}, we have proved the following uniqueness result:

\begin{theorem}\label{unique}
The dimension of the space of left integrals on $\A$ is at most equal to
\nolinebreak 1.
In particular, every integral on $\A$ is even or odd.
\end{theorem}

As in the classical case of the Haar measure on locally compact groups,
this theorem enables us to determine how left integrals behave under
``right translations".

\begin{proposition}\label{rint}
Let $\int$ be a non-trivial left integral on a Hopf superalgebra $\A$\,.
Then there exists a unique even group-like element $a_0 \in \A$ such that
 \[ (\tint \ot\, \id_{\A}) \comp \Delta = a_0 \tint \;. \]
In particular, $\int$ is also a right integral if and only if $a_0 = 1_\A$\,.
\end{proposition}

For later reference, we mention the following simple lemma.

\begin{lemma}\label{auxres}
(a) Let $f : \A \rar {\cal B}$ be an {\em injective} homomorphism of a Hopf
superalgebra $\A$ into a Hopf superalgebra $\cal B$. If $\int$ is a left
integral on $\cal B$, then $\int \comp f$ is a left integral on $\A$\,. \\
(b) In particular, if $\A$ is a sub-Hopf-superalgebra of a Hopf superalgebra
$\cal B$, then the restriction of a left integral on $\cal B$ onto $\A$ is
a left integral on $\A$\,. \\
(c) Let $\A_1,\ldots,\A_n$ be $n$ Hopf superalgebras. Then the left integrals
on the tensor product $\A_1 \ot \ldots \ot \A_n$ (endowed with the tensor
product Hopf superalgebra structure) are exactly the linear forms of the type
$\int_1 \ot \ldots \ot \int_n$\,, where, for every $i$\,, $\int_i$ is a left
integral on $\A_i$\,.
\end{lemma}

\renewcommand{\thesection}{\arabic{section}.\hspace{-0.5em}}
\sect{Some comments on the finite dual of a Hopf superalgebra \vspace{-1ex}}
\renewcommand{\thesection}{\arabic{section}}
In the present work we are mainly interested in the case where $\A$ is a
sub-Hopf-superalgebra of the finite (or continuous) dual $\Uc$ of a Hopf
superalgebra $\U$. We recall that, as a vector space, $\Uc$ consists of all
linear forms on $\U$ that vanish on a graded two-sided ideal of $\U$ of
finite codimension. It is easy to see that the properties derived in
Ref.~\cite{Swe} can immediately be extended to the present setting. In
particular, this is true for the criteria characterizing the elements of
$\Uc$. Moreover, the elements of $\Uc$ can also be described by the following
property: \\[1.0ex]
A linear form $\ell$ on $\U$ belongs to $\Uc$ if and only if there exists a
finite-dimensional graded (left) $\U$-module $V$, an element $v \in V$, and
a linear form $v^{\ast} \in V^{\ast}$ such that
\be
  \ell(x) = \la v^{\ast},\tty x \cd v \ra \quad\mbox{for all $x \in \U$} \,,
                                                         \label{findual} \ee
where the dot denotes the module action, and where $\la \;,\: \ra$ denotes
the canonical pairing of $V^{\ast}$ and $V$. By a slight abuse of language,
the linear forms $\ell$ of the type \reff{findual} are called the {\em matrix
elements} of the graded $\U$-module $V$ (or of the associated representation
of $\U$ in $V$).

Obviously, $\Uc$ is a graded subspace of $\U^{\ast}$, the graded dual of $\U$.
Actually, much more is true, for the structure of $\U$ as a Hopf superalgebra
yields a similar structure on $\Uc$. More precisely, the transpose of
$\De_{\U}$ induces the multiplication in $\Uc$, the counit $\ve_{\U}$ is the
unit element of $\Uc$, the transpose of the product map of $\U$ induces the
coproduct of $\Uc$, the evaluation at the unit element $1_{\U}$ of $\U$ is
the counit of $\Uc$, and the transpose of $S_{\U}$ induces the antipode of
$\Uc$. For simplicity, the comultiplication, counit, and antipode of $\U$
will also be denoted by $\Delta\tty$, $\ve$, and $S$, respectively.

Suppose now that $\U$ is a Hopf superalgebra, and that $\cO$ is a class of
finite-dimensional graded $\U$-modules that contains at least one non-trivial
module and is closed under the formation of finite direct sums, tensor
products and duals. Then the set of all matrix elements of modules belonging
to $\cO$ is a sub-Hopf-superalgebra $\A_{\cO}$ of $\Uc$. If $V$ belongs to
$\cO$, there exists a unique structure of a graded right $\A_{\cO}$-comodule
on $V$ such that the structure map
 \[ \de : V \lra V \ot \A_{\cO} \]
satisfies the following condition: If $v$ is an element of $V$ and if
 \[ \de(v) = \sum_i v^0_i \ot v^1_i \;, \]
with homogeneous elements $v^0_i \in V$ and $v^1_i \in \A_{\cO}$\tty, then
for every homogeneous element $x \in \U$, the module action of $x$ on $v$ is
given by
\be   x \cd v \,=\, (-1)^{[x][v]} \sum_i v^0_i \la v^1_i , x \ra \;. 
                                                         \label{comod} \ee

Conversely, let $\A$ be a sub-Hopf-superalgebra of $\Uc$, and let $\cO'$ be
the class of all finite-dimensional graded right $\A$-comodules. If $V$
belongs to $\cO'$, Eqn.~\reff{comod} converts the graded vector space $V$
into a graded $\U$-module $V^{\rm mod}$. Denoting the class of all graded
$\U$-modules $V^{\rm mod}$ obtained in this way by $\cO$, it can be shown
that $\cO$ is closed under the formation of finite direct sums, tensor
products and duals, and that the corresponding sub-Hopf-superalgebra
$\A_{\cO}$ is equal to $\A$\,. (To prove the last statement, regard $\A$ as
a graded right $\A$-comodule via $\De$\,, and then use the fact that every
element of $\A$ is contained in a finite-dimensional graded subcomodule
\linebreak
of $\A$\,.)

Keeping the notation of the preceeding paragraph, it is easy to see that a
subset $W$ of $V$ is a graded $\A$-subcomodule of $V$ if and only if it is
a graded $\U$-submodule of $V^{\rm mod}$. In particular, this implies that a
finite-dimensional graded right $\A$-comodule $V$ is simple or semisimple if
and only if the graded $\U$-module $V^{\rm mod}$ has this property.

Occasionally, it is adequate to assume that $\A$ is {\em dense} in
$\U^{\ast}$. By definition, this is to say that for every non-zero element
$x \in \U$ there exists an element $a \in \A$ such that $\la a,x \ra \neq 0$
(using yet another language, this means that $\A$ separates the elements of
$\U$\,). This assumption appears to be quite natural. Indeed, let $K$ be the
intersection of the kernels of all elements of $\A$\,. It is not difficult to
see that $K$ is a graded Hopf ideal of $\U$. Hence if $K$ is different from
$\Ze$, we are essentially considering the quotient Hopf superalgebra $\U/K$
rather than $\U$ itself.

\renewcommand{\thesection}{\arabic{section}.\hspace{-0.5em}}
\sect{A method to construct integrals on Lie \\
      supergroups \vspace{-1ex}}
\renewcommand{\thesection}{\arabic{section}}
In the present section we are going to recall a technique, developed in
Ref.~\cite{SZh}, by which integrals on Lie supergroups can be constructed.

Let $\fg = \fgz \op \fgo$ be a finite-dimensional Lie superalgebra
\cite{Kac, Sch}. We take $\U$ to be the enveloping algebra $U(\fg)$ of
$\fg$\,. $U(\fg)$ contains the enveloping algebra $U(\fgz)$ of the Lie
subalgebra $\fgz$ as a sub-Hopf-algebra. We denote $U(\fgz)$ by $\Ue$\,, and
let
 \[ \cI:\, \Ue \lra \U \]
be the canonical embedding, which is a Hopf superalgebra map. It is well-known
that the transpose $\cI^{\ast}$ of $\cI$ induces a Hopf superalgebra map
 \[ \cP : \Uc \lra \Uec \;, \vspace{-1.5ex} \]
which is given by
 \[ \la\ty \cP(a), u \ra = \la a\,,\,\cI(u)\ra
                      \quad\mbox{for all $a \in \Uc$, $u \in \Ue$} \;. \]

In the present work, a Lie supergroup will be defined in terms of its Hopf
superalgebra of functions, i.e., we proceed as in the usual definition of
quantum groups \cite{RTF} or quantum supergroups \cite{Man} (for a related
treatment of supergroups, see \cite{Kos} and \cite{Bos}). More precisely,
if $\fg$ is a Lie superalgebra, the superalgebra of functions on a Lie
supergroup associated to $\fg$ will be a sub-Hopf-superalgebra $\A$ of
$\Uc = U(\fg)^{\circ}$, subject to the condition that $\A$ be dense in
$U(\fg)^{\ast}$. Actually, in our discussion of integrals, this latter
property will be used for motivation purposes only.

Thus, let $\A$ be a sub-Hopf-superalgebra of $\Uc$. We set
 \[ \cP(\A) =  \Ae \;, \] 
which is a Hopf subalgebra of $\Uec$. Moreover, let
 \[ \nu : \U = U(\fg) \lra  \A^{\ast} \]
be the canonical map, given by
 \be \la \nu(x)\,,\,a \ra \,=\, (-1)^{[x][a]} \la a \,,\, x \ra
                                                          \;, \label{nu} \ee
for all homogeneous elements $x \in \U$ and $a \in \A$\,. Actually, $\nu$
maps $\U$ into $\Ac$, and the induced map of $\U$ into $\Ac$ is a Hopf
superalgebra map (which is injective if $\A$ is dense in $\U^{\ast}$).

Next, we define
 \[ J = \U \fgz \;. \]
Obviously, $J$ is a graded left ideal of $\U$, and the Poincar\'e, Birkhoff,
Witt theorem implies that its codimension in $\U$ is finite. Consequently,
$\U/J$ is a finite-dimensional graded left $\U$-module. We note that this
module is isomorphic to the $\U$-module induced from the trivial $\Ue$-module.

Using these data, we have the following theorem, proved in Ref.~\cite{SZh}.

\begin{theorem}\label{main}
Let $\int_e$ be a left integral on $\Ae$\tty, and let $z + J \in \U/J$ be
invariant under the action of $\U$. Then we have: \\
(a) The linear form $\int = \nu(z) \cd (\int_e \comp \cP)$ is a left integral
on $\A$ and does not depend on the choice of the representative for $z + J$. \\
(b) Suppose in addition that $\int_e 1_{\Ae} \!\neq 0$, that the matrix
elements of the $\U$-module $\U/J$ belong to $\A$\,, and that $z \notin J$.
Then the integral $\int$ of part (a) is different from zero.
\end{theorem}

We recall that the dot in part (a) of the theorem denotes the product in
$\A^{\ast}$ deduced from the coproduct in $\A$\,. It is worthwhile to notice
that
 \[ \tint\, 1_\A \,=\, \ve(z) \tinte\, 1_{\Ae} \,. \]
Consequently, $\int 1_\A$ is different from zero if and only if $\ve(z)$
and $\int_e 1_{\Ae}$ are different from zero.

In order to find a non-trivial integral on $\A$ by application of
Theorem \ref{main}, we have to satisfy the conditions of part (b). The
inequality $\int_e 1_{\Ae} \!\neq  0$ is very restrictive. According to a
classical theorem (see Ref.~\cite{Swe}), this is equivalent to the requirement
that all $\Ae$-comodules are semisimple. \vspace{2ex}

\noindent
{\em Remark 4.1.} The aforementioned result is a special case of a much
more general theorem by Larson \cite{Lar}. In particular, Larson's theorem
implies that an analogous result holds for Hopf superalgebras as well. For
these, the theorem has been rederived in Ref.~\cite{SZh}. %\vspace{1ex}

\noindent
{\em Remark 4.2.} At this point we should also mention the following result
by Sullivan \cite{Sul}: If $H$ is a {\em commutative} Hopf algebra, and if
$\int$ is a non-zero integral on $H$, then it follows that
$\int 1_H \neq 0\ty$ (provided the characteristic of the base field is equal
to zero). Thus there is no point in trying to relax the condition
$\int_e 1_{\Ae} \!\neq  0$ by just assuming that $\int_e$ is different from
zero. \vspace{2ex}

As mentioned at the beginning of this section, we are finally interested in
the case where $\A$ is dense in $\U^{\ast}$, and hence $\Ae$ is dense in
$\Ue^{\ast}$. According to the discussion in Section 3 (applied to $\Ae$ and
$\Ue$), this implies that $\Ue$ admits a class of finite-dimensional
semisimple modules, whose matrix elements separate the elements of $\Ue$\tty.
In particular, this shows that $\fgz$ has a faithful finite-dimensional
completely reducible representation. According to a well-known theorem in Lie
algebra theory, the latter statement characterizes the reductive Lie algebras.

In view of this result, {\em we assume until the end of the present section
that $\fgz$ is a reductive Lie algebra,} i.e., that
 \[ \fgz = \fs \op \fz \quad\mbox{as Lie algebras} \;, \]
where $\fs = [\fgz\tty,\fgz]$ is semisimple and $\fz$ is the center of
$\fgz\tty$. It follows that
 \[ U(\fgz) \simeq U(\fs) \ot U(\fz) \quad\mbox{as Hopf algebras} \]
and hence that
 \be U(\fgz)^{\circ} \simeq U(\fs)^{\circ} \ot U(\fz)^{\circ}
                          \quad\mbox{as Hopf algebras} \;. \label{Uciso} \ee
Now we intend to apply Lemma \ref{auxres}(c). It is well-known that on
$U(\fs)^{\circ}$ there exists a left (and right) integral $\int_{\fs}$ such
that $\int_{\fs} 1_{U(\fs)^{\circ}} = 1\,$. (Essentially, this follows from
the theorem on cosemisimple Hopf algebras mentioned above; see also the
discussion in Ref.~\cite{SZh}.) On the other hand, if $\fz \neq \Ze$, there
do not exist any non-trivial integrals on $U(\fz)^{\circ}$ (use Lemma
\ref{auxres}(c) and the Appendix of Ref.~\cite{SZh}). However, let $G$ be the
set of grouplike elements of $U(\fz)^{\circ}$, i.e., the set of characters of
the algebra $U(\fz)$. Under the multiplication in $U(\fz)^{\circ}$, this is a
group, isomorphic to the underlying additive group $\fz_{+}^{\ast}$ of the
dual $\fz^{\ast}$ of $\fz\tty$. Hence the linear span $\CC\tty G$ of $G$ is a
sub-Hopf-algebra of $U(\fz)^{\circ}$, and isomorphic to the group Hopf algebra
of $\fz_{+}^{\ast}$\tty. Thus there exists a left (and right) integral
$\int_{\fz}$ on $\CC\tty G$ such that $\int_{\fz} 1_{\CC\tty G} = 1\,$.

Let $\B$ be the sub-Hopf-algebra of $U(\fgz)^{\circ}$ that under the
isomorphism \reff{Uciso} corresponds to $U(\fs)^{\circ} \ot \CC\tty G$. Then
the foregoing discussion shows that there exists a left (and right) integral
$\int_{\B}$ on $\B$ such that
 \[ \textstyle{\int_{\B} 1_{\B} = 1 \;.} \]
It is useful to observe that $\B$ has a natural interpretation in terms of
$U(\fgz)$-modules: Using the notation of Section 3, $\B$ is the
sub-Hopf-algebra $\A_{\cS}$ of $U(\fgz)^{\circ}$ associated to the class
$\cS$ of all finite-dimensional semisimple $U(\fgz)$-modules, i.e., $\B$
consists exactly of the matrix elements of the finite-dimensional semisimple
$U(\fgz)$-modules.

Now let $\cO$ be the class of all finite-dimensional graded $U(\fg)$-modules
which are semisimple when regarded as $U(\fgz)$-modules. Obviously, $\cO$
is closed under the formation of finite direct sums, tensor products and
duals. In view of the preceding result, we choose $\A$ to be the
sub-Hopf-superalgebra $\A_{\cO}$ of $U(\fg)^{\circ}$ as defined in Section 3,
i.e., the sub-Hopf-superalgebra consisting of the matrix elements of those
finite-dimensional graded $U(\fg)$-modules which are semisimple when regarded
as $U(\fgz)$-modules. By definition, $\Ae$ is contained in $\B$. Consequently,
there exists a left (and right) integral $\int_e$ on $\Ae$ such that
 \[ \tinte\, 1_{\Ae} = 1 \,, \] 
and hence the first condition on $\A$ is satisfied.

The second condition to be satisfied is that the matrix elements of the
$\U$-mod\-ule $\U/J$ should belong to $\A$\,. Obviously, this is the case if
the $\Ue$-module $\U/J$ is semisimple, and this is true if (and only if)
the $\fgz$-module $\fgo$ is semisimple. Moreover, if the $\fgz$-module $\fgo$
is semisimple, it also follows that $\Ae = \B$, and that $\A$ is dense in
$\U^{\ast}$ (use the theory of induced modules \cite{Sch} and argue
similarly as in the proof of Theorem 2.5.7 in Ref. \cite{Dix}).

Summarizing, we have shown that with our choice for $\A$\,, the first and the
second condition in Theorem \ref{main}(b) are satisfied, provided that the
Lie algebra $\fgz$ is reductive and that the $\fgz$-module $\fgo$ is
semisimple. It remains to investigate under which conditions the graded
$\U$-module $\U/J $ contains non-zero invariants. As we are going to see in
the subsequent section, this problem can be solved quite generally.

\renewcommand{\thesection}{\arabic{section}.\hspace{-0.5em}}
\sect{Existence of integrals \vspace{-1ex}}
\renewcommand{\thesection}{\arabic{section}}
The following lemma contains the solution to the problem stated at the end
of the preceding section.

\begin{lemma}\label{basic}
Let $\fg = \fgz \op \fgo$ be a finite-dimensional Lie superalgebra, let
$\U$ be its enveloping algebra, and let ${\rm ad}'$ be the representation
of $\fgz$ in $\fgo$ given by the adjoint representation of $\fg\ty$. Then
the subspace of invariant elements of the $\U$-module $\U/\U\fgz$ is at most
one-dimensional, and it is one-dimensional if and only if
 \be {\rm Tr}({\rm ad}'(X)) = 0 \quad\mbox{for all $X \in \fgz$} \;.
                                                       \label{tracecond} \ee
\end{lemma}

{\em Proof\,}: Let $\lam$ be the linear form on $\fgz$ defined by
 \[ \lam(X) = {\rm Tr}({\rm ad}'(X)) \quad\mbox{for all $X \in \fgz$} \;. \]
Using the universal property of the enveloping algebra $\Ue$ of $\fgz\tty$,
it is easy to see that there exists a unique automorphism $\al$ of the
algebra $\Ue$ such that
 \[ \al(X) = X + \lam(X)\ty 1_{\Ue} \quad\mbox{for all $X \in \fgz$} \;. \]

Now let $W$ be the vector space consisting of all linear forms $\ell$ on
$\U$ such that
 \[ \ell(u_e u) = \ve(\al(u_e)) \ell(u)
                    \quad\mbox{for all $u_e \in \Ue$ and $u \in \U$} \,, \]
where $\ve$ denotes the counit of $\U$. We make $W$ into a left $\U$-module
by means of the action
 \[ (u \cd \ell)(u') = \ell(u' u)
                  \quad\mbox{for all $\ell \in W$ and $u,u' \in \U$} \,. \]
(Note that we have not introduced a $\ZZ_2\tty$-gradation on $W$, although we
could easily do so.) Then the theory of Frobenius extensions \cite{NTs, BFa} 
implies that the (non-graded) $\U$-modules $\U/\U\fgz$ and $W$ are
isomorphic (see the Appendix for more details). Consequently, it is
sufficient to determine the $\U$-invariant elements of $W$.

This is easily done. By definition, an element $\ell \in W$ is invariant
if and only if
 \[ u \cd \ell = \ve(u)\ty \ell \quad\mbox{for all $u \in \U$} \,, \]
i.e., if and only if
 \[ \ell(u' u) = \ell(u')\tty \ve(u) \quad\mbox{for all $u,u' \in \U$} \,. \]
Obviously, this is true if and only if $\ell$ is proportional to $\ve\tty$.
On the other hand, $\ve$ belongs to $W$ if and only if
 \[ \ve(u_e) = \ve(\al(u_e)) \quad\mbox{for all $u_e \in U_e$} \;, \]
i.e., if and only if
 \[ \lam(X) = 0 \quad\mbox{for all $X \in \fgz$} \;. \]
This proves the lemma. In  Eqn.~\reff{final} of the appendix, we give an
explicit (albeit rather complicated) formula for the invariant $z + \U\fgz$
corresponding to $\ve\tty$ (assuming that the condition \reff{tracecond} is
satisfied). \vspace{2.0ex}

Summarizing our results, we have proved the following theorem.

\begin{theorem}\label{exist}
Let $\fg$ be a finite-dimensional Lie superalgebra. Suppose that the Lie
algebra $\fgz$ is reductive, and that the linear map ${\rm ad}'(Z)$ is
diagonalizable and traceless, for all central elements $Z$ of $\fgz\tty$.
Define the sub-Hopf-superalgebra $\A$ of $U(\fg)^{\circ}$ as at the end of
Section 4. Then $\A$ is dense in $U(\fg)^{\ast}$, and there exists a non-zero
left integral on $\A$\,. The integral is unique up to scalar multiples.
\end{theorem}

{\em Proof\,}: All we have to do is to recall that (for a reductive Lie
algebra $\fgz$) a finite-dimensional $\fgz$-module is semisimple if and only
if the representatives of the central elements of $\fgz$ are diagonalizable.
\vspace{2.0ex}

\noindent
{\bf Corollary}\hspace{0.5em}
{\em Let $\fg$ be a classical simple Lie superalgebra. Define the
sub-Hopf-superalgebra $\A$ of $U(\fg)^{\circ}$ as at the end of Section 4.
Then $\A$ is dense in $U(\fg)^{\ast}$, and there exists a non-zero left
integral on $\A$\,. The integral is unique up to scalar multiples, and it is
also a right integral.} \vspace{2.0ex}

{\em Proof\,}: The assumptions of Theorem \ref{exist} are satisfied. Hence
it remains to show that the left integral is also a right integral. Since the
Lie superalgebra $\fg$ is simple, the counit $\ve$ is the sole superalgebra
homomorphism of $\U$ into $\CC\ty$. According to Proposition \ref{rint}, this
proves our claim. (See also Lemma 2 of Ref.~\cite{SZh}.)

\renewcommand{\thesection}{\arabic{section}.\hspace{-0.5em}}
\sect{Discussion \vspace{-1ex}}
\renewcommand{\thesection}{\arabic{section}}
In the present work, we have proved the existence of non-trivial integrals
for a reasonably large class of Lie supergroups (or, rather, on their Hopf
superalgebras of functions). A look at the defining formula shows that,
qualitatively speaking, the integrals are defined by ``first integrating over
the even parameters of the supergroup, and then over the odd ones". The
reasoning of Section 4 seems to indicate that with this type of approach one
cannot treat much more general cases than we did.

The integration over the odd parameters is governed by the invariant $z + J$.
In the appendix, we shall give a formula for $z + J$ in terms of a so-called
dual free pair relative to an associative form. Unfortunately, this formula
seems to be rather complicated, and it remains to be seen whether it is
useful in practical calculations. In any case, we would like to understand
the structure of $z + J$ much better than we do at present. It is worth
noting that this invariant plays a decisive role in Gorelik's theory of the
so-called anticenter of $U(\fg)$ \cite{Gor}.  \vspace{4.0ex}

\noindent
{\LARGE\bf Appendix\vspace{-1.5ex}}

\begin{appendix}

\renewcommand{\thesection}{\Alph{section}.\hspace{-0.5em}}
\section{Definitions and results from the theory of \\
         Frobenius extensions \vspace{-1ex}}
\renewcommand{\thesection}{\Alph{section}}
The theory of Frobenius extensions \cite{NTs, BFa} has been an important tool
in the proof of Lemma \ref{basic} and hence of Theorem \ref{exist}. Therefore,
we think it may be useful to give a short summary of those definitions and
results of the theory that are directly relevant to our application. We shall
follow the presentation given in Ref.~\cite{BFa} (indeed, for the convenience
of the reader we shall even use the notation of that reference). Moreover,
we shall elaborate a little on some points which have been mentioned in
Ref.~\cite{BFa} very briefly, but which are crucial here.

Throughout this appendix all rings are assumed to have a unit element that
operates via the identity transformation on all modules under consideration.
Let $R$ be a ring, and let $S$ be a subring of $R$ containing the unit
element of $S$. Furthermore, let $\al$ be an automorphism of the ring $S$.
If $M$ is a left $S$-module, we let ${}_{\al}M$ denote the left $S$-module
which, regarded as an abelian group, is equal to $M$, but endowed with the
new action of $S$ defined by
 \[ s \ast m = \al(s)m \quad\mbox{for all $s \in S$ and $m \in M$} \,. \]
For any two left $S$-modules $M$ and $M'$, the abelian group (under addition)
of all $S$-module homomorphisms of $M$ into $M'$ will be denoted by
$\Hom_S^{\ell}(M,M')$ (while for right modules we use the superscript $r$).
In particular, $\Hom_S^{\ell}(R\tty,{}_{\al}S)$ is the set of additive maps
$f: R \rar S$ such that
 \[ f(sr) = \al(s)f(r) \quad\mbox{for all $s \in S$ and $r \in R$} \;. \]
This is an $(R\tty,S)$-bimodule via the action
 \[ (r \cd f \cd s)(x) \,=\, f(xr)s 
                       \quad\mbox{for all $r,x \in R$ and $s \in S$} \,. \]
Note that an additive map $f : R \rar S$ lies in $\Hom_S^{\ell}(R\tty,S)$ if
and only if $\al \comp f$ lies in $\Hom_S^{\ell}(R\tty,{}_{\al}S)$. We say
that $R$ is a {\em free $\al$-Frobenius extension} of $S$ if
\begin{enumerate}
\item[(${\rm i}_{\ell}$)]
 $R$ is a finitely generated free left $S$-module, and \vspace{-0.5ex}
\item[(${\rm ii}_{\ell}$)]
 there exists an isomorphism
 $\vp_{\ell}: R \rar \Hom_S^{\ell}(R\tty,{}_{\al}S)$ of $(R\tty,S)$-bimodules.
\end{enumerate}

The $\al$-Frobenius extensions have a nice description in terms of
$\al^{-1}$-associative forms from $R$ to $S$, in the sense of the following
definition. Let $\bet$ be a ring endomorphism of $S$.
A {\em $\bet$-associative form from $R$ to $S$} is a biadditive map
$\la \;,\: \ra : R \ti R \rar S$ such that
 \[ \la sx , y \ra = s \ty \la x , y \ra \quad,\quad
            \la xr , y \ra = \la x , ry \ra \quad,\quad
                    \la x , ys \ra = \la x , y \ra \bet(s) \]
for all $s \in S$ and $x,y,r \in R$\tty. We say that $\la \;,\: \ra$ is
{\em non-degenerate} if its left and right radical are equal to $\Ze$.

A $\bet$-associative form from $R$ to $S$ determines an additive map
$\pi : R \rar S$ by $\pi(r) = \la r , 1 \ra = \la 1 , r \ra$ that satisfies
 \be \pi(sr) = s \, \pi(r) \quad\mbox{and}\quad \pi(rs) = \pi(r)\bet(s)
     \quad\mbox{for all $r \in R$ and $s \in S$} \, . \label{Frobhomid} \ee
Conversely, every additive map $\pi : R \rar S$ with these properties
yields a $\bet$-associa\-tive form from $R$ to $S$ by means of
$\la x,y \ra = \pi(xy)$ for all $x,y \in R$\tty.

Suppose now that $R$ is an $\al$-Frobenius extension of $S$, and let
$\vp_{\ell}$ be the isomorphism mentioned in condition (${\rm ii}_{\ell}$).
We define a map $\la \;,\: \ra : R \ti R \rar S$ by
 \[ \la x,y \ra = \al^{-1}(\vp_{\ell}(y)(x))
                                     \quad\mbox{for all $x,y \in R$} \;, \]
or, eqivalently, by
 \[ \vp_{\ell}(y)(x) = \al(\la x,y \ra) \quad\mbox{for all $x,y \in R$} \;. \]
Then it is easy to check that $\la \;,\: \ra$ is an $\al^{-1}$-associative
form from $R$ to $S$.

Let $\{x_1,\ldots,x_n\}$ be a basis of the left $S$-module $R\tty$ (the
existence is guaranteed by condition (${\rm i}_{\ell}$)), and let
$f_i \in \Hom_S^{\ell}(R\tty,S)$, $1 \leq i \leq n$\tty, be the associated
coordinate forms, defined by the condition that
 \[ r \,= \sum_{i=1}^n f_i(r)\ty x_i \quad\mbox{for all $r \in R$} \;. \]
As noted above, we have $\al \comp f_i \in  \Hom_S^{\ell}(R\tty,{}_{\al}S)$. 
Let $y_i \in R$ be the element such that $\vp_{\ell}(y_i) = \al \comp f_i$\ty.
Using the fact that $\vp_{\ell}$ is an isomorphism, it can be shown that
$\{x_1,\ldots,x_n\}$ and $\{y_1,\ldots,y_n\}$ form a dual free pair relative
to $\la \;,\: \ra$, in the sense of the following definition.

Let $\la \;,\: \ra : R \ti R \rar S$ be an $\al^{-1}$-associative form. Two
subsets $\{x_1,\ldots,x_n\}$ and $\{y_1,\ldots,y_n\}$ of $R$ are said to form
a {\em dual free pair relative to} $\la \;,\: \ra$ if
\begin{enumerate}
\item[(i)]  $R = \sum_{i=1}^n S x_i = \sum_{i=1}^n y_i \tty S$,
\item[(ii)] $\la x_i\ty,y_j \ra = \de_{ij}$ for $1 \leq i\tty, j \leq n$\ty.
\end{enumerate}
Obviously, this implies that $\{x_1,\ldots,x_n\}$ and $\{y_1,\ldots,y_n\}$
are bases of $R$ regarded as a left and right $S$-module, respectively, and
that the form $\la \;,\: \ra$ is non-degenerate.

As shown in Ref.~\cite{BFa}, the converse of the above is also true, i.e.,
we have the following theorem.

\begin{theorem}\label{Frob}
Let $S$ be a subring of $R$ (containing the unit element of $R$) and let
$\al$ be an automorphism of the ring $S$. Then the following statements are
equivalent: \\
(a) $R$ is a free $\al$-Frobenius extension of $S$. \\
(b) There is an $\al^{-1}$-associative form $\la \;,\: \ra$ from $R$ to $S$
relative to which a dual free pair $\{x_1,\ldots,x_n\}$, $\{y_1,\ldots,y_n\}$
exists.
\end{theorem}

The map $\pi : R \rar S$ defined by $\pi(r) = \la r,1 \ra = \la 1,r \ra$
will be called the {\em Frobenius homomorphism} of the Frobenius extension.
Recall that the form $\la \;,\: \ra$ is determined by $\pi\tty$.

Now let $R^{\rm\ty op}$ and $S^{\rm op}$ be the rings opposite to $R$ and $S$,
respectively. Then the criterion (b) of Theorem \ref{Frob} shows that $R$
is a free $\al$-Frobenius extension of $S$ with form $\la \;,\: \ra$ if and
only if $R^{\rm\ty op}$ is a free $\al^{-1}$-Frobenius extension of
$S^{\rm op}$ with the form $\la \;,\: \ra_{\rm op}$ defined by
 \[ \la x,y \ra_{\rm op} = \al(\la y,x \ra)
                                      \quad\mbox{for all $x,y \in R$} \;. \]

We want to reformulate the fact that $R^{\rm\ty op}$ is a free
$\al^{-1}$-Frobenius of $S^{\rm op}$ in terms of $R$ and $S$ themselves. To
do so we need some definitions. If $M$ is a right $S$-module, we let
$M_{\al^{-1}}$ denote the right $S$-module which, regarded as an abelian
group, is equal to $M$, but endowed with the new action of $S$ defined by
 \[ m \ast s = m \,\al^{-1}(s)
                          \quad\mbox{for all $m \in M$ and $s \in S$} \,. \]
For any two right $S$-modules $M$ and $M'$, the abelian group (under addition)
of all $S$-module homomorphisms of $M$ into $M'$ will be denoted by
$\Hom_S^r(M,M')$. In particular, $\Hom_S^r(R\tty,S_{\al^{-1}})$ is the set of
additive maps $f: R \rar S$ such that
 \[ f(rs) = f(r)\al^{-1}(s) \quad\mbox{for all $r \in R$ and $s \in S$} \;. \]
This is an $(S,R)$-bimodule via the action
 \[ (s \cd f \cd r)(x) \,=\, sf(rx) 
                       \quad\mbox{for all $s \in S$ and $r,x \in R$} \,. \]

Using these definitions, the reformulation we are looking for reads as
follows: Let $R$ be a ring and let $S$ be a subring of $R$ (containing the
unit element of $R$). Then $R$ is a free $\al$-Frobenius extension of $S$ if
and only if
\begin{enumerate}
\item[(${\rm i}_r$)]
 $R$ is a finitely generated free right $S$-module, and \vspace{-0.5ex}
\item[(${\rm ii}_r$)]
 there exists an isomorphism
 $\vp_r : R \rar \Hom_S^r(R\tty,S_{\al^{-1}})$ of $(S,R)$-bimodules.
\end{enumerate}
The equivalence of the system of conditions (${\rm i}_{\ell}$),
(${\rm ii}_{\ell}$) with the system (${\rm i}_r$), (${\rm ii}_r$) has been
mentioned in Remark (1) on page 412 of Ref.~\cite{BFa}, and it has also been
proved in Ref.~\cite{NTs}. It is the second system that will be used in the
following. Let us note that the isomorphisms $\vp_{\ell}$ and $\vp_r$ can be
chosen such that
 \[ \vp_r(x)(y) = \al^{-1}(\vp_{\ell}(y)(x)) = \la x,y \ra 
                                     \quad\mbox{for all $x,y \in R$} \;. \]

Now let $M$ be a left $S$-module, and let $R \ot_S M$ be the corresponding
induced left $R$-module. Quite generally, there exists a canonical
$R$-module homomorphism
 \be \tau : R \ot_S M
        \lra \Hom_S^{\ell}(\Hom_S^r(R\tty,S_{\al^{-1}}),{}_{\al}M) \;,
                                                             \label{tau} \ee
defined by
 \be (\tau(r \ot_S m))(f)
                       \,=\, ((\al \comp f)(r))\ty m \;, \label{tauexpl} \ee
for all $r \in R\tty$, $m \in M$, $f \in \Hom_S^r(R\tty,S_{\al^{-1}})$. (Some
comments are in order. The left action of $R$ on the right hand side of
\reff{tau} is derived, in the usual way, from the right action of $R$ on
$\Hom_S^r(R\tty,S_{\al^{-1}})$. On the right hand side of
Eqn.~\reff{tauexpl}, the element $(\al \comp f)(r) \in S$ acts on $m \in M$
in the sense of the $S$-module $M$ (not ${}_{\al}M$).) If condition
(${\rm i}_r$) is satisfied, the homomorphism $\tau$ is bijective, i.e., an
$R$-module isomorphism. But then condition (${\rm ii}_r$) implies that
there also exists an $R$-module isomorphism
 \be \om : R \ot_S M \lra \Hom_S^{\ell}(R\tty,{}_{\al}M) \;. \label{iso} \ee
This result is due to Ref.~\cite{NTs}. Recalling the definitions of the
various maps involved, we find that $\om$ is given by
 \be (\om(r \ot_S m))(r') \,=\, \al(\la r',r \ra)\ty m
       \quad\mbox{for all $r,r' \in R$ and $m \in M$} \,. \label{omexpl} \ee
Note that for $M = S$ we recover the isomorphism $\vp_{\ell}\tty$.

In the present work, we are mainly interested in the following special case.
Suppose that $R$ is an algebra (associative, with unit element) over a field
$\KK$\,, and that $S$ is a subalgebra of $R$ containing the unit element of
$R\tty$. Let $\ve$ be a character of $R$\tty, i.e., an algebra homomorphism
of $R$ into $\KK$ (mapping $1_R$ onto $1_{\KK}$). We regard $\KK$ as a left
$S$-module by means of the action
 \[ s \cd c = \ve(s)\tty c \quad\mbox{for all $s \in S$ and $c \in \KK$} \]
and want to determine all elements $\ol{z}$ of the induced $R$-module
$R \ot_S \KK$ such that
 \be r \cd \ol{z} = \ve(r)\ty \ol{z}
                         \quad\mbox{for all $r \in R$} \;. \label{cond} \ee

Using the isomorphism \reff{iso}, we first determine the elements
$f \in \Hom_S^{\ell}(R\tty,{}_{\al}\KK)$ such that
 \[ r \cd f = \ve(r)\ty f \quad\mbox{for all $r \in R$} \;. \]
Proceeding exactly as in the proof of Lemma \ref{basic}, we see that non-zero
elements $f$ of this type exist if and only if
 \be \ve(\al(s)) = \ve(s) \quad\mbox{for all $s \in S$} \;. \label{inv} \ee
Moreover, if this condition is satisfied, the elements in question are the
linear forms proportional to $\ve\tty$.

Suppose now that \reff{inv} is true. We have to determine the preimage
$\om^{-1}(\ve)$ of $\ve$ in $R \ot_S \KK$\,. Let $\{x_1,\ldots,x_n\}$,
$\{y_1,\ldots,y_n\}$ be a dual free pair relative to $\la \;,\: \ra$. Since
$\{y_1,\ldots,y_n\}$ is a basis of the right $S$-module $R$\tty, every
element $v \in R \ot_S \KK$ can uniquely be written in the form
 \[ v \,= \sum_{i=1}^n y_i \ot_S c_i \;, \]
where $c_i \in \KK$ for $1 \leq i \leq n$\tty. According to
Eqn.~\reff{omexpl}, we have for all elements $r \in R$
 \[ \begin{array}{rcll}
   \dis{(\om(\ty\sum_i y_i \ot_S c_i))(r)}
%       &\!=\!& \dis{\sum_i
%          \la r,y_i \ra \cd c_i} & \mbox{action in ${}_{\al}\KK$}  \\[2.5ex]
       &\!=\!& \dis{\sum_i
          \al(\la r,y_i \ra) \cd c_i} & \mbox{action in $\KK$}     \\[2.5ex]
       &\!=\!& \dis{\sum_i
          \ve(\al(\la r,y_i \ra))\tty c_i} &                       \\[2.5ex]
       &\!=\!& \dis{\sum_i \ve(\la r,y_i \ra)\tty c_i \;.} &
    \end{array} \]
If $v$ is the preimage of $\ve\tty$, this has to be equal to $\ve(r)$, for
all $r \in R\tty$. Evaluating on $r = x_j\ty$, we obtain $c_j = \ve(x_j)$,
for $1 \leq j \leq n\tty$. Consequently, we have
 \be \om^{-1}(\ve) \,= \sum_{i=1}^n y_i \ot_S \ve(x_i) \;, \label{invar} \ee
and the elements $\ol{z} \in R \ot_S \KK$ satisfying Eqn.~\reff{cond} are
exactly those proportional to $\om^{-1}(\ve)$. 
 
It remains to make contact with the theory of Lie superalgebras. The following
theorem is due to Bell and Farnsteiner \cite{BFa}.

\begin{theorem}
Let $\KK$ be a field whose characteristic is different from $2\tty$, let
$\fg = \fgz \op \fgo$ be a finite-dimensional Lie superalgebra, and let
$\al$ be the algebra automorphism of $U(\fgz)$ defined in the proof of
Lemma \ref{basic}. Then $U(\fg)$ is a free $\al$-Frobenius extension of
$U(\fgz)$.
\end{theorem}

Actually, Bell and Farnsteiner have proved a more general result, but the
special case given by the theorem above is sufficient for our purposes.
Actually, they first prove this special case by constructing the Frobenius
homomorphism $\pi$ (hence the $\al^{-1}$-associative form $\la \;,\: \ra$)
and a dual free pair relative to this form. It is worthwhile to recall at
least part of these constructions.

Let $(x_i)_{1 \leq i \leq m}$ be a basis of $\fgo\tty$, and let
${\cal J} = \{0,1\}^m$ be the set of all $m$-tuples with entries in
$\{0,1\}$. In particular, we put $T = (1,\ldots,1)$. For each element
$I = (i(1),\ldots,i(m)) \in {\cal J}$, we define
${\bf x}^I = x_1^{i(1)} \!\ldots\ty x_m^{i(m)} \in U(\fg)$. Then
$({\bf x}^I)_{I \in {\cal J}}$ is a basis both of the left and right
$U(\fgz)$-module $U(\fg)$. Accordingly, every element $u \in U(\fg)$ can be
written uniquely in the form
 \[ u = \sum_{I \in {\cal J}} u_I \tty {\bf x}^I \;,
         \quad\mbox{with $u_I \in U(\fgz)$ for all $I \in {\cal J}$} \,. \]
We define
 \be \pi : U(\fg) \lra U(\fgz) \label{Frobhom} \ee
by
 \[ \pi(\sum_{I \in {\cal J}} u_I \tty {\bf x}^I) \,=\, u_T \;. \]  
It is easy to see that $\pi$ satisfies the conditions \reff{Frobhomid} with
$\bet = \al^{-1}$. Hence the associated form $\la \;,\: \ra$ is
$\al^{-1}$-associative. Let us also note that if we choose another basis
of $\fgo\tty$, the map $\pi$ and hence the form $\la \;,\: \ra$ will change
by a non-zero scalar factor.

For every element $I \in {\cal J}$, let $I' = T - I$ denote the
``complement" of $I$, and set $|I| = \sum_{r=1}^m i(r)\ty$. Choose a total
order $\ll$ on $\cal J$ such that $I \ll J$ implies $|I| \leq |J|\tty$.
Then the matrix
 \[ A \,=\, (\la {\bf x}^I , {\bf x}^{J'} \ra)_{I,\ty J \in {\cal J}} \]
is lower triangular, and its entries on the diagonal are all equal to
$\pm 1\ty$. This implies that $A$ is invertible. Defining
 \[ {\bf y}^I \,= \sum_{J \in {\cal J}} {\bf x}^{J'} \al((A^{-1})_{JI})
                               \quad\mbox{for all $I \in {\cal J}$} \,, \]
we conclude that $\{{\bf x}^I | I \in {\cal J}\}$ and 
$\{{\bf y}^I | I \in {\cal J}\}$ form a dual free pair relative to the form
$\la \;,\: \ra$.

Let us now assume that Eqn.~\reff{inv}, i.e., Eqn.~\reff{tracecond} is
satisfied, and let $N = T' = (0,\ldots,0)$ be the smallest element of
$\cal J$. Observing that $\ve({\bf x}^I) = 0$ if $I \neq N$, we conclude
from Eqn.~\reff{invar} that the invariant $\om^{-1}(\ve)$ is given by
 \be \om^{-1}(\ve) \,=\, {\bf y}^N \ot_S 1_{\KK}
     \,= \, \sum_{J \in {\cal J}} {\bf x}^{J'} \ot_S \ve((A^{-1})_{JN}) \;.
                                                           \label{final} \ee

In the appendix we have not yet considered the $\ZZ_2\tty$-gradations.
Obviously, the Frobenius homomorphism \reff{Frobhom} is even or odd depending
on whether $m = \dim \fgo$ is even or odd. It follows that the same holds
true for the $\al^{-1}$-associative form $\la \;,\: \ra$ and for the module
isomorphism $\om\ty$. Thus $\om$ is an isomorphism of graded $\U(\fg)$-modules
if $\dim \fgo$ is even, but a shift of the gradations is involved if
$\dim \fgo$ is odd.

\end{appendix}

\end{document}